\newcommand\webcite[1]{\texttt{\def~{\~{}}#1}}
\documentclass{article}
\usepackage{amsmath,amsfonts,amsthm}
\newtheorem{theorem}{Theorem}
\newtheorem{lemma}[theorem]{Lemma}

\renewcommand\Pr{{\mathop{\mathbb P{}}\nolimits}}
\def\Z{{\mathbb Z}}

\def\gs{X}
\def\eps{\varepsilon}
\newcommand\downto{\searrow}

\begin{document}
\title{A note on the Harris-Kesten Theorem}
\date{September 10, 2005}

\author{B\'ela Bollob\'as\thanks{Department of Mathematical Sciences,
University of Memphis, Memphis TN 38152, USA}
\thanks{Trinity College, Cambridge CB2 1TQ, UK}
\thanks{Research supported in part by NSF grant ITR 0225610}
\and Oliver Riordan%
\thanks{Royal Society Research Fellow, Department of Pure Mathematics
and Mathematical Statistics, University of Cambridge, UK}}
\maketitle

\begin{abstract}
A short proof of the Harris-Kesten result that the critical
probability for bond percolation in the planar square lattice is $1/2$
was given in \cite{ourKesten}, using a sharp threshold result of
Friedgut and Kalai. Here we point out that a key part of this proof
may be replaced by an argument of Russo~\cite{Russo01} from 1982,
using his approximate zero-one law in place of the Friedgut-Kalai result.
Russo's paper gave a new proof of the Harris-Kesten Theorem
that seems to have received little attention.
\end{abstract}

Let $\Z^2$ be the planar square lattice, i.e., the graph with vertex
set $\Z^2$ in which each pair of nearest neighbours is joined by an
edge. Let $\gs=E(\Z^2)$ be the edge-set of $\Z^2$, and let
$\Omega=\{-1,+1\}^\gs$. We write $\omega=(\omega_e)_{e\in \gs}$ for an
element of $\Omega$, and say that the edge $e$ is {\em open} (in the
state $\omega$) if $\omega_e=+1$, and {\em closed} if $\omega_e=-1$.
An event $A\subset \Omega$ is {\em local} if it depends on
only finitely many coordinates. As usual, let $\Sigma$ be the sigma-field
generated by local events, and let $\Pr_p$ be the probability measure
on $(\Omega,\Sigma)$ in which each edge is open with probability $p$,
and these events are independent.
Let $\theta(p)$ be the $\Pr_p$-probability that the origin is in an
{\em infinite open cluster},
i.e., an infinite connected subgraph $C$ of $\Z^2$ with every edge
of $C$ open.
In 1960, Harris~\cite{Harris} proved that $\theta(1/2)=0$; in 1980,
Kesten~\cite{Kesten1/2} showed that $\theta(p)>0$ for $p>1/2$, establishing
that $p_c=1/2$ is the `critical probability' for this model.
A short proof of these results was given in \cite{ourKesten},
using a sharp-threshold result of Friedgut and Kalai~\cite{FK}, itself
based on a result of Kahn, Kalai and Linial~\cite{KKL}.

In 1982, Russo~\cite{Russo01} proved a general sharp-threshold result
(weaker than the more recent results described above)
and applied it to percolation,
to give a new proof of the `equality of critical probabilities' for
site percolation in $\Z^2$. Although Russo does not explicitly say this,
his application applies equally well to bond percolation, giving
a new proof of the Harris-Kesten Theorem that seems not to be well known.
Here we shall present Russo's general sharp-threshold result,
and then give a complete version of his application,
to bond percolation in $\Z^2$.

Replacing the appropriate section of \cite{ourKesten} with this argument
gives an even simpler proof of the Harris-Kesten Theorem; we are grateful
to Professor Ronald Meester for bringing this to our attention.

An event $A\subset \Omega$ is {\em increasing} if $\omega\in A$ and
$\omega_e\le \omega_e'$ for
every $e$ imply $\omega'\in A$, i.e., if
$A$ is preserved when the state of one or more edges is changed from closed to open.
An edge $e$ is {\em pivotal} for an event $A$ if changing the state of $e$ affects
whether or not $A$ holds. Let $\delta_e A$ be the event
that $e$ is pivotal for $A$, so $\omega\in \delta_e A$ if and only if exactly one
of $\omega^+,\omega^-$ is in $A$, where $\omega^{\pm}$ are the states that agree
with $\omega$ on all edges other than $e$, with $\omega^+_e=1$ and $\omega^-_e=-1$.
In \cite{Russo01}, Russo proved the following result about the product measure $\Pr_p$;
in this result the structure of $\Z^2$ is irrelevant, i.e.,
the groundset $X$ can be any countable set.

\begin{theorem}\label{th_R}
For every $\eps>0$ there is an $\eta>0$ such that if $A$ is an increasing local event
with
\[
 \Pr_p(\delta_e A)<\eta
\]
for every $e\in \gs$ and every $p\in [0,1]$, then there is a $p_0\in [0,1]$ with
\[
 \Pr_{p_0-\eps}(A) \le \eps \hbox{ and } \Pr_{p_0+\eps}(A)\ge 1-\eps.
\]
\end{theorem}

As in \cite{ourKesten},
by a {\em $k$ by $\ell$ rectangle} we mean
a rectangle $[a,b]\times[c,d]$ with $a,b,c,d\in \Z$ and $b-a=k$, $d-c=\ell$.
We identify a rectangle with the corresponding subgraph of $\Z^2$, including the boundary.
A rectangle $R$ has a {\em horizontal open crossing} if there is a path in $R$
consisting of open edges, joining a vertex on the left-hand side of $R$ to one on the right;
we write $H(R)$ for this event. Our starting point will be the following
consequence of the Russo-Seymour-Welsh Lemma (see \cite{ourKesten} and the references
therein): there is a constant $c>0$ such that
\begin{equation}\label{start}
 \Pr_{1/2}(H(R)) \ge c,
\end{equation}
for any $3n$ by $n$ rectangle $R$.
This is essentially the case $\rho=3$ of Corollary 7 in \cite{ourKesten}. (The latter
result has an irrelevant restriction to $n$ even; the present statement is immediate from
the case $\rho=4$ of this result.)

Our aim is to deduce Lemma 11 of \cite{ourKesten}, restated below.
\begin{lemma}\label{l_end}
Let $p>1/2$ be fixed.
If $R_n$ is a $3n$ by $n$ rectangle, then
$ \Pr_p(H(R_n)) \to 1$
as $n\to \infty$.
\end{lemma}
It is well known that Lemma \ref{l_end} implies Kesten's Theorem; see~\cite{ourKesten}.
We shall deduce Lemma \ref{l_end} from \eqref{start} using
Theorem \ref{th_R} and Harris' result, that $\theta(1/2)=0$. We shall
need the concept of the {\em dual lattice} $(\Z^2)^*$: this is the planar dual of the graph
$\Z^2$, having a vertex for each face of $\Z^2$, and an edge $e^*$ for each edge
$e$ of $\Z^2$, joining the two vertices corresponding to the faces of $\Z^2$ in whose boundary
$e$ lies. We take $e^*$ to be open if and only if $e$ is closed. The following argument
is based on that of Russo~\cite{Russo01}.

\begin{proof}[Proof of Lemma \ref{l_end}]
Let $p_1>1/2$ be fixed.
Let $D$ be a constant to be chosen below, and let $R$ be a $3n$ by $n$ rectangle with
$n\ge 2D+1$.
Suppose that $\omega\in \delta_e H(R)$, and define $\omega^{\pm}$
as above. Note that $e$ must be an edge of $R$, as $H(R)$ depends only on such edges.
Then, in $\omega^+$ there is an open path in $R$ from the left-hand side to the right
using the
edge $e$. Hence, in $\omega$, the endpoints of $e$ are joined by open paths
to the left- and right-hand sides of $R$. One of these paths must have length at least
$(3n-1)/2\ge D$.
Thus, for any $p$,
\begin{equation}\label{e1}
 \Pr_p(\delta_e H(R)) \le 2\Pr_p(0\to D),
\end{equation}
where $0\to D$ is the event that there is an open path of length $D$ starting at the origin.
Our assumption that $e$ is pivotal also implies that $H(R)$ does not hold in $\omega^-$.
It follows (by Lemma 3 of \cite{ourKesten}) that in $\omega^-$ there is an open path in the
dual
lattice joining the top of $R$ to the bottom, using the edge $e^*$. Hence, in the dual
lattice, one of the endpoints of $e^*$ is in an open path of length at least $D$.
As edges of the dual lattice are open independently with probability $1-p$, it follows that
\begin{equation}\label{e2}
 \Pr_p(\delta_e H(R)) \le 2\Pr_{1-p}(0\to D).
\end{equation}

Let $0<\eps<\min\{(p_1-1/2)/2,c\}$ be arbitrary, where $c>0$ is a constant for which
\eqref{start}
holds.
Let $\eta=\eta(\eps)$ be as in Theorem \ref{th_R}.
For any $p$ we have $\Pr_p(0\to D)\downto \theta(p)$ as $D\to \infty$.
Hence, by Harris' Theorem (Theorem 8 in \cite{ourKesten}), $\Pr_{1/2}(0\to D)\to 0$,
so we may choose $D$ such that $\Pr_{1/2}(0\to D)\le \eta/3$.
As the event $0\to D$ is increasing, for $p\le 1/2$ we have
\[
 \Pr_p(0\to D)\le \Pr_{1/2}(0\to D)\le \eta/3.
\]
Using \eqref{e1} for $p\le 1/2$ and \eqref{e2} for $p\ge 1/2$, it follows that for any
$p\in [0,1]$ and any edge $e$ in $R$ we have
\[
 \Pr_p(\delta_e H(R))\le 2\eta/3<\eta.
\]
As $H(R)$ is an increasing local event, and $\delta_e H(R)$ is empty for edges outside
$R$, the conditions of Theorem \ref{th_R} are satisfied.
Hence, $\Pr_p(H(R))$ increases from at most $\eps<c$ to at least $1-\eps$ in some interval
of width at most $2\eps<p_1-1/2$. As $\Pr_{1/2}(H(R))\ge c$ by \eqref{start},
it follows that $\Pr_{p_1}(H(R))\ge 1-\eps$.
In other words, we have shown that for $p_1>1/2$ and $\eps>0$ fixed and $R_n$ a $3n$ by $n$
rectangle,
we have $\Pr_{p_1}(H(R_n))\ge 1-\eps$ if $n$ is large enough.
As $\eps>0$ is arbitrary, this completes the proof.
\end{proof}

In Section 5 of \cite{ourKesten}, the Friedgut-Kalai sharp threshold
result is used to deduce from \eqref{start} a result (Lemma 9 in
\cite{ourKesten}) that is somewhat stronger than Lemma \ref{l_end}.
This stronger form was used in the first proof of Kesten's Theorem
given in \cite{ourKesten}; however, in \cite{ourKesten} two more very
simple proofs are given, both of which need only Lemma \ref{l_end}.

\medskip
\noindent{\bf Acknowledgement.}
We would like to thank Professor Ronald Meester for drawing Russo's paper
to our attention, and pointing out that Russo's proof may replace the
relevant argument in \cite{ourKesten}.

\end{document}